\documentclass[11pt,twoside]{article}
\usepackage[utf8]{inputenc}
\usepackage{epsfig,amsfonts,color}
\usepackage{mathrsfs} 
\usepackage{amsmath, bbm, mathabx, amsthm}
\usepackage[normalem]{ulem}
\usepackage{amssymb, palatino, geometry,url}
\usepackage[colorlinks=true,linkcolor=blue,citecolor=blue,urlcolor=blue]{hyperref}
\usepackage{subcaption}
\usepackage{enumitem}
\usepackage{multirow}

\usepackage[title]{appendix}
\usepackage{subfiles}

\usepackage{algorithm}
\usepackage{algpseudocode}
\usepackage{array}
\usepackage{mathtools}

\usepackage{threeparttable} 
\newcounter{subsubsubsection}[subsubsection]

\usepackage{booktabs}
\usepackage{longtable}
\usepackage{fancyhdr}
\usepackage{lineno}
\usepackage{float}
\usepackage[section]{placeins}
\usepackage[table]{xcolor}

\title{Optimal Production Scheduling for \\Multi-Product Mineral Manufacturing}
\author{Stefan Mazzadi$^1$ and Joshua L. Pulsipher$^{1}$\thanks{Corresponding Author: pulsipher@uwaterloo.ca}\\
    {\small $^1$Department of Chemical Engineering}\\
    {\small \;University of Waterloo, 200 University Ave W, Waterloo, ON N2L 3G1, Canada}
    }
\date{}

\begin{document}
\maketitle
\begin{abstract}
Short-term scheduling in multi-product mineral manufacturing requires coordinating material movement and production across interconnected processing stages while addressing confirmed customer orders and longer-term production requirements. This work develops a generalized network-based, lexicographic multi-objective optimization (MOO) scheduling framework for these systems. The formulation represents the progression of material from raw-material supply through production and inventory to customer fulfillment. A sequential objective hierarchy first prioritizes fulfillment of confirmed customer orders and then progressively directs feasible production capacity toward cumulative targets supplied by the facility's demand-planning process. This approach ensures that the performance achieved for each higher-priority objective is preserved before subsequent targets are considered. Selected product-line activation decisions are represented using generalized disjunctive programming (GDP). An industrial mineral-manufacturing case study demonstrates the proposed framework's ability to generate feasible shift-level schedules and preserve the intended operational priorities with practical computational performance.
\end{abstract}

\noindent\textbf{Keywords:} Production Scheduling; Multi-Objective Optimization; Generalized Disjunctive Programming; Mineral Manufacturing

\section{Introduction}
\label{sec:intro}
Production scheduling is a critical activity in large-scale manufacturing where multiple product grades, interacting equipment units, material inventories and operational restrictions must be coordinated to meet demand efficiently \cite{Gaglioppa_2008,mendez2006state}. These challenges become particularly significant in multi-product processing facilities where product-dependent production rates, product changeovers, storage limitations, equipment compatibility and shared resources can strongly influence the feasibility and performance of a schedule. Mineral processing operations exhibit many of these characteristics, with materials of different grades and physical properties moving through interconnected storage and production resources before reaching finished-product inventories. Developing effective schedules for such systems requires simultaneous consideration of material availability, production capacity, product changeover requirements, inventory limitations and competing longer-term operational priorities.

Mathematical programming provides a systematic approach for representing these interconnected scheduling decisions and constraints. M\'{e}ndez et al. \cite{mendez2006state} provide an extensive review of mixed integer formulations for short-term scheduling of batch processes, highlighting the importance of appropriately representing material flows, equipment assignments and sequencing decisions. They also explore the trade-offs between discrete and continuous time formulations, and computational complexity. 
More recently, Samadi and Maravelias \cite{samadi2024comprehensive} introduce the general material task system representation for describing multiple production environments, material-handling restrictions, transient operations, and storage policies. Optimization-based scheduling has similarly been applied to multi-stage manufacturing systems characterized by large product portfolios, interacting processing resources, and packaging operations \cite{georgiadis2018production}. For industrial systems in which operating conditions and production requirements evolve over time, online scheduling approaches provide a means of repeatedly updating scheduling decisions as new plant information becomes available \cite{gupta2019online}. Despite these methodological advances, translating optimization-based scheduling methods into industrial decision-support systems remains challenging, as ensuring deployability, computational tractability, maintainability, and compatibility with existing information systems is imperative \cite{harjunkoski2014scope,stuve2022systematic}.

Production scheduling can also involve multiple competing performance criteria that cannot always be represented adequately through a single objective. Depending on the production system, relevant criteria may include customer service, tardiness, inventory, resource utilization, production cost, workload distribution, and the frequency or duration of product changeovers. Multi-objective optimization (MOO) provides several approaches for representing competing criteria, including weighted aggregation, Pareto-based methods and lexicographic optimization \cite{Ehrgott2005, Miettinen1998}. Multi-objective formulations have been applied directly to production scheduling problems. T'kindt et al. \cite{TKINDT200142} formulate glass production scheduling as a bicriteria problem involving production margin and workload balance, while Sawik \cite{SAWIK20071060} applies a lexicographic approach to bi-objective scheduling in make-to-order manufacturing. Almada-Lobo et al. \cite{ALMADALOBO2008363} consider multiple scheduling criteria associated with sequence-dependent setup times, inventory levels and stockouts in glass container production, while Demir \cite{DEMIR2023110541} applies lexicographic MOO to an industrial fabric-dyeing scheduling problem in which tardiness, sequence-dependent washing requirements and machine cost are prioritized according to operational preferences. More recently, de Sousa \cite{desousa2026development} applies a two-stage lexicographic mixed integer linear program (MILP) to multi-plant glass container production planning, prioritizing demand fulfillment before transport cost. Collectively, these studies demonstrate the use of different multi-objective strategies for reconciling competing production priorities. Lexicographic MOO is particularly suitable when objectives possess an explicit priority order, since higher-order priority objectives can be optimized first and their attained performance preserved as subsequent objectives are considered \cite{Ehrgott2005, Miettinen1998, DEMIR2023110541}.

Generalized disjunctive programming (GDP) provides a structured mathematical framework for representing such logical relationships through Boolean propositions and algebraic disjunctions \cite{perez2023julia}. GDP has also been applied to selected scheduling problems in which logical and discrete operating decisions play an important role. Castro and Grossmann \cite{CastroGrossmann2012} use GDP as a systematic framework for deriving short-term batch process scheduling formulations, while Castro and Marques \cite{castro2015operating} apply continuous-time GDP formulations to operating-room scheduling. More recently, Wu et al. \cite{WU2021107268} develop a GDP-based formulation for joint production and maintenance scheduling in a multi-product batch chemical plant, incorporating sequence-dependent degradation and limited finished product storage. Software tools such as \texttt{DisjunctiveProgramming.jl} \cite{perez2023julia} and \texttt{Pyomo.GDP} \cite{chen2022pyomo} readily facilitate the translation of GDP constraints to forms that are compatible with mixed-integer program (MIP) solvers. In the present work, GDP is employed as a modelling mechanism for representing logical scheduling decisions rather than as the principal methodological contribution.

Within the mining industry, optimization-based scheduling has been applied to extraction and mine-to-plant material-flow decisions, including mine sequencing, loading and haulage allocation, stockpile management, ore blending and the delivery of material to processing plants \cite{MARTINS2024110115}. Liu et al. \cite{liu2025mining} develop an integrated open-pit production-timetabling framework that optimizes the sizes of mining jobs with the allocation and sequencing of equipment across multiple operational stages. Recent work has also incorporated multiple operational targets at increasingly detailed time scales. Martins et al. \cite{MARTINS2024110115} integrate a hierarchical multi-objective MILP formulation with discrete-event simulation to generate shift-by-shift short-term plans for an open-pit iron-ore operation, prioritizing production, ore-grade, and particle size targets for processing plants together with mine-level operating objectives. These approaches demonstrate the value of optimization for coordinating extraction, material preparation, and mine-to-plant material flows. However, the studies considered above focus primarily on upstream mining operations rather than the detailed downstream scheduling of multi-product mineral processing and manufacturing networks. Downstream operations require coordination of different mineral grades through interconnected storage and production resources while accounting for resource compatibility, product changeovers, shared equipment, finished-product inventories, customer demand and longer-term production targets.

To address this gap, this work develops a general multi-objective production scheduling optimization framework for multi-product mineral-processing systems. The proposed framework represents the flow of material from source nodes through intermediate storage, production resources and finished-product inventories while incorporating product-resource compatibility, production capacities, product changeovers, shared equipment interactions, inventory limitations and customer-demand requirements. Discrete operating decisions are represented using GDP logical constraints. A lexicographic MOO strategy establishes an explicit hierarchy among operational priorities pertaining to customer-demand satisfaction and production targets over successive time horizons. The principal contributions of this work include the following:
\begin{itemize}
    \item the development of a generalized network-based representation for multi-product mineral-processing scheduling,
    \item the use of a hierarchical lexicographic MOO strategy for prioritizing customer-demand satisfaction and sequential product-target attainment,
    \item and the application of the above contributions to a sector-specific industrial-scale mineral manufacturing case study.
\end{itemize}
The remainder of this paper is structured as follows. Section \ref{sec:problem_definition} defines the mineral-processing scheduling problem. Section \ref{sec:methodology} outlines the proposed scheduling formulation and solution strategy. Section \ref{sec:casestudy} applies the framework to an industrial case study. Section \ref{sec:conclusion} provides concluding remarks.

\section{Problem Definition} \label{sec:problem_definition}
This work considers the short-term scheduling of a multi-product mineral manufacturing facility in which bulk raw material passes through intermediate storage and production resources before being stored as finished product and shipped to customers. The scheduling problem determines how material should be routed through this network, which products should be produced on each production line, how much of each product should be produced, where the resulting production should be stored, and how available finished-product inventory should be allocated to customer orders throughout the horizon. These decisions are strongly interdependent as production lines differ in product compatibility and processing rate, production lines compete for shared equipment, material can move only via a set of shared transfer paths, and both raw material and finished-product storage capacities are limited. In addition to maintaining operational feasibility, the schedule must prioritize confirmed customer orders and longer-term production requirements.

\subsection{Scheduling Horizon and Time Representation}\label{subsec:time_representation}
The scheduling horizon is divided into an ordered set of discrete periods $t\in\mathcal{T}$, referred to as production shifts. The scheduling problem is organized around production shifts as staffing, operating records, and production commitments are managed at shift boundaries. The duration of a shift can be selected according to the operating practices specific to the facility. Material transfers, production campaigns, and customer shipments are assigned to shifts, while inventory levels are observed at the corresponding boundaries. This resolution supports coordination over a multi-day horizon without requiring the exact start time or within-shift sequence of each activity. The required output is therefore a shift-level operating plan rather than a detailed equipment-level execution schedule.

\subsection{Manufacturing Network and Resource Connectivity}\label{subsec:manufacturing_network}
\begin{figure}[!htb]
	\centering
	\includegraphics[width=1\textwidth]{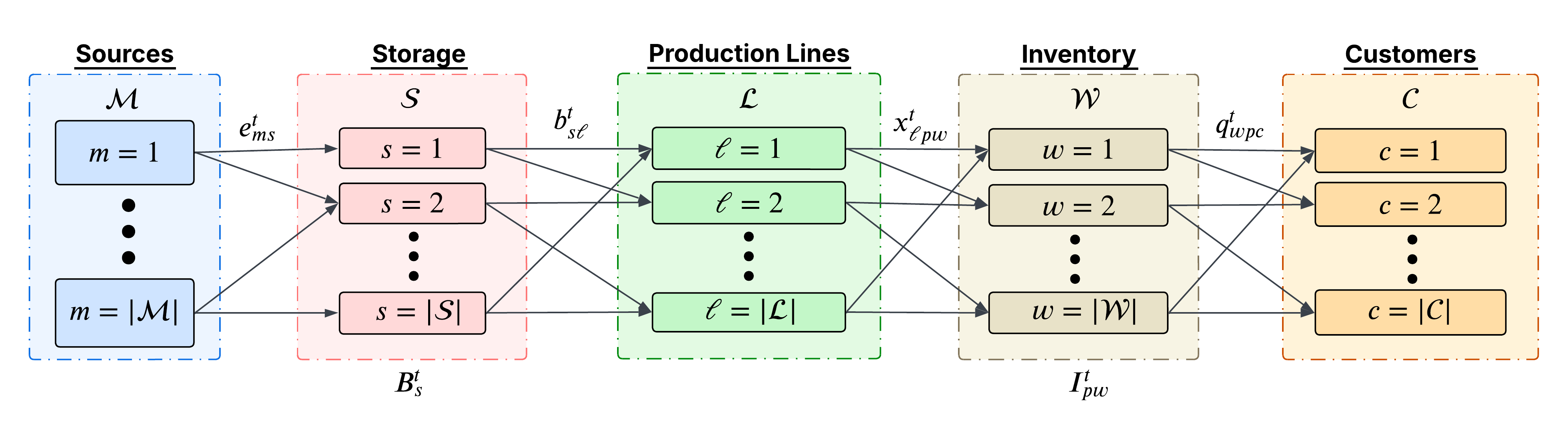}
	\caption{Illustration of the proposed manufacturing network abstraction.}
	\label{fig:PFD}
\end{figure}

The facility is represented as a five-layer material-flow network comprising raw-material sources, raw-material storage nodes, production lines, finished-product inventory locations and customers, as illustrated in Figure \ref{fig:PFD}. These layers are represented by the sets $\mathcal{M}$, $\mathcal{S}$, $\mathcal{L}$, $\mathcal{W}$, and $\mathcal{C}$, respectively. Finished products are indexed by $p\in\mathcal{P}$. This network abstraction retains the physical connectivity, compatibility relationships and capacity limitations that determine whether a schedule is operationally feasible without explicitly modelling the detailed mechanical operation of every conveyor, hopper, processing unit or packaging component involved in the process.

Each source $m\in\mathcal{M}$ represents an aggregate supply point from which material produced by an upstream extraction or processing route becomes available to the manufacturing facility. The upstream production operations themselves are outside the scope of the scheduling problem, where the source nodes define the points at which material from these upstream operations becomes available to the downstream manufacturing network.

Material available from different sources may differ in purity, particle size, or handling characteristics, which affect where material can be stored and which finished products it can support. The subset $\mathcal{M}_s\subseteq\mathcal{M}$ identifies the raw-material source nodes capable of supplying material to raw-material storage node $s$, while the reciprocal subset $\mathcal{S}_m\subseteq\mathcal{S}$ contains the raw-material storage nodes capable of receiving material from source node $m$.

The set $\mathcal{S}$ represents bulk silos, holding bins and other raw-material holding locations. These storage nodes decouple upstream mineral availability from downstream production but are subject to finite capacities and grade-segregation requirements. Storage nodes and production lines are not fully interconnected as restricting piping, conveying systems, elevation differences and equipment layouts limit the available transfer routes. Additionally, material characteristics impose additional restrictions as a storage node containing one mineral grade may not be suitable for a line producing a different grade or product specification. For example, a fine-grain, high-purity material may be routed only through segregated storage and production equipment qualified for that grade. The subset $\mathcal{S}_{\ell}\subseteq\mathcal{S}$ contains the storage nodes capable of supplying line $\ell$, while the reciprocal subset $\mathcal{L}_s\subseteq\mathcal{L}$ contains the lines capable of receiving material from storage node $s$. Material can therefore be routed only through connections supported by the facility layout and compatible with the relevant mineral grade and product requirements.

The production line set $\mathcal{L}$ represents the processing and packaging resources used to transform bulk mineral into finished-product formats. Products in $\mathcal{P}$ may differ in mineral grade, purity, particle size, treatment requirements, package format, and customer specifications. Consequently, neither the products nor the production lines can be treated as interchangeable resources. The subset $\mathcal{P}_{\ell}\subseteq\mathcal{P}$ contains the products compatible with line $\ell$, while $\mathcal{L}_p\subseteq\mathcal{L}$ contains the lines capable of producing product $p$. These compatibility relationships reflect equipment capabilities, material-handling requirements, packaging configurations, product-quality requirements, and achievable production rates. Even when a product can be manufactured on more than one line, its achievable production rate may differ between lines because of differences in equipment capacity and product-specific processing requirements.

Production-line availability may also be coupled through equipment that supports more than one nominally distinct line. For example, two lines may depend on a common upstream hopper, conveyor, screening unit or shared downstream material-handling equipment. Production assigned to one line can reduce the operating time available to another even though the lines are represented as separate production resources. The set $\mathcal{L}_{\mathrm{conflicts}}\subseteq\mathcal{L}\times\mathcal{L}$ contains the unordered pairs of lines subject to such shared-resource interactions. The shared resource is required only while a connected line is actively producing, as line-specific changeover activities do not occupy the shared equipment and may occur while another connected line continues production.

The operating time available during each shift is reduced by planned downtime, operator shift changeovers and other non-productive activities. Product-dependent changeovers consume additional time within this operating window, meaning the shift-level schedule must also account for the operating implications of changing production between products. Each activated product corresponds to a production campaign that may require line preparation, equipment configuration, packaging changes, cleaning, purging or other transition activities. 

After processing, finished product is assigned to an inventory location $w\in\mathcal{W}$. These locations may represent pre-customer holding bins, quality-hold staging area, specific warehouse zones, or other facility packaged-product inventory locations. The subset $\mathcal{W}_p\subseteq\mathcal{W}$ identifies the locations capable of storing product $p$, while $\mathcal{P}_{w}\subseteq\mathcal{P}$ identifies the products compatible with inventory location $w$. These relationships capture restrictions associated with product package type, storage configuration, contamination prevention and established site-specific inventory practices. Each location is subject to aggregate lower and upper inventory bounds that can represent physical capacity, safety-stock requirements and other operating policies. Together, these reciprocal subsets describe the same physical compatibility relationships from either side of each network connection.

\subsection{Customer Demand and Forecast-Based Production Priorities} \label{subsec:scheduling_priorities}
Finished-product inventory is used to satisfy confirmed customer orders. Customers are indexed by $c\in\mathcal{C}$ and the set $\mathcal{O}\subseteq\mathcal{C}\times\mathcal{P}$ contains the customer--product pairs for which demand exists within the scheduling horizon. Demand outstanding before the scheduling horizon constitutes an initial backorder, while demand arising during the horizon becomes eligible for fulfillment on its corresponding due shift. Customer-service performance therefore depends on both the quantity of confirmed demand that remains unmet and the duration for which it remains outstanding.

In addition to confirmed customer orders, the facility must work toward forecast-based production targets intended to support anticipated future requirements. In practice, these product-level targets may be established through demand-planning processes that generate statistical forecasts from historical demand data and incorporate judgmental adjustments reflecting anticipated customer or market expectations \cite{FILDES20093}. Although the future requirements underlying these targets are uncertain in nature, the resulting target quantities are treated as fixed planning inputs for each scheduling horizon. Development of the production forecast and the explicit representation of its uncertainty are outside the scope of the scheduling problem. When the scheduling horizon begins partway through an existing forecast cycle, production already completed within that cycle is credited toward the corresponding forecast targets. The subset $\mathcal{P}^{\mathrm{forecast}}\subseteq\mathcal{P}$ identifies the products included in the production forecast, and the ordered set $\mathcal{K}=\{1,\ldots, N\}$ represents the successive forecast intervals. 

Confirmed customer demand and forecast-based production targets each possess explicitly different operational priorities. Demand that has already become due takes precedence over forecast attainment, after which remaining production capacity is directed toward forecast intervals in chronological order. Production supporting an earlier forecast interval should not be displaced solely to improve attainment of a later interval.

The required schedule specifies shift-level material transfers, product-line assignments, production quantities, finished-product inventory allocations, and customer shipments. Upstream extraction, detailed within-shift sequencing, individual material-handling movements, operator assignments, and scheduling of maintenance activities are outside the problem scope. Demand, forecast targets, initial inventories, processing rates, equipment compatibility, and available capacities are assumed known over the scheduling horizon. Material quantities are treated as continuous, processing losses are assumed negligible, and unplanned equipment failures or staffing disruptions are not considered.

\section{Model Formulation and Solution Methodology} \label{sec:methodology}
Building on the manufacturing-network representation and scheduling requirements defined in Section \ref{sec:problem_definition}, this section introduces the corresponding mathematical formulation. The formulation first describes the material balances, production capacities, logical activation decisions, inventory restrictions and customer demand constraints governing a feasible shift-level schedule. The resulting scheduling model is then embedded within a sequential lexicographic MOO procedure that preserves customer-fulfillment performance while progressively considering forecast production targets. The principal sets, indices, parameters, and decision variables used throughout the formulation are summarized in Section \ref{sec:notation}.

\subsection{Material Balance and Capacity Constraints}
The formulation follows the five network layers illustrated in Figure \ref{fig:PFD}, proceeding from raw-material sources through intermediate storage and production lines to finished-product inventory and customer fulfillment. Flow variables represent quantities transferred, produced, or shipped during a shift, while inventory variables represent the corresponding material levels accumulated at the designated shift boundaries. The following subsections define the variables and constraints associated with each network layer.

\subsubsection{Raw Material Source Nodes}\label{subsubsec:source_constraints}
Raw-material source nodes $m\in\mathcal{M}$ represent the upstream points at which raw mineral enters the modelled production network. As mentioned in Section \ref{subsec:manufacturing_network}, each source in a mineral-manufacturing system may correspond to a mine or another upstream extraction or processing route that provides a specified mineral feed to the facility. As the upstream extraction and processing operations are outside the scope of the scheduling model, the non-negative parameter $E_m^t$ specifies the total quantity of raw material supplied by source $m$ during shift $t$. This quantity is provided as an exogenous model input and may be obtained from a fixed delivery schedule, a short-term mine production plan or the expected generation rate of the corresponding upstream operations. Consequently, the scheduling model does not determine how much raw material is produced at the source.

The non-negative continuous variable $e_{ms}^t$ denotes the portion of the supply transferred to compatible raw material storage node $s\in\mathcal{S}_m$ during shift $t$. Although raw material may physically arrive through individual deliveries or production batches, the model represents material as an aggregate mass quantity over each shift and does not schedule individual trucks, conveyance events, or batches. As a result, the continuous representation is appropriate when the supplied mineral can be treated as divisible in the resolution of the shift-level scheduling. Since the supplied quantity is predetermined and must be allocated entirely among the compatible storage nodes, the following source balance is imposed:
\begin{equation}
\sum_{s\in\mathcal{S}_m} e_{ms}^t = E_m^t
\quad
\forall m\in\mathcal{M},\;t\in\mathcal{T}
\label{eq:source_availability}
\end{equation}
Equation \eqref{eq:source_availability} enforces material conservation at the source layer by requiring the aggregate source-to-storage flows to equal the supplied quantity. The optimization model determines only how the fixed supply from each source is distributed among its compatible storage nodes. These source-to-storage transfers provide the inflows to the raw-material storage balances defined in the following subsection.

\subsubsection{Raw Material Storage Nodes}\label{subsubsec:storage_constraints}
Raw-material storage nodes $s\in\mathcal{S}$ represent the silos, holding bins, and other bulk-storage units positioned between the upstream mineral sources and the production lines. These nodes buffer the fixed supply entering from the source layer and make accumulated material available to downstream production over the subsequent shifts. Each storage node is associated with mineral characteristics such as grade, purity, or particle size, and may receive material only from sources for which $s\in\mathcal{S}_m$. Similarly, material may be withdrawn from storage node $s$ only by production lines $\ell\in\mathcal{L}_s$ that are physically connected to the node and compatible with its mineral contents.

In this model, it is assumed that a source-to-storage transfer initiated during shift $(t-1)$ requires a single shift to reach its destination and consequently becomes available for storage and production during shift $t$. For each storage node $s$, a material balance accounts for the inventory carried forward from the preceding shift, the raw material received from compatible sources during the preceding shift, and the material withdrawn by compatible production lines:
\begin{equation}
  B_s^t = B_s^{t-1} + \sum_{m \in \mathcal{M}_s}e_{ms}^{t-1} - \sum_{\ell \in \mathcal{L}_s} b^t_{s\ell} \quad \forall s \in \mathcal{S}, t \in \mathcal{T}
  \label{eq:storage_balance}
\end{equation}
where $B^t_{s} \in [0,B_s^{\max}]$ is the inventory level of raw material at storage node $s$ at the end of shift $t$, $e_{ms}^{t-1}$ is the inflow from raw material source $m$ to storage node $s$ during shift $t-1$, and the non-negative continuous variable $b^t_{s\ell}$ is the flow from storage node $s$ to production line $\ell$ during shift $t$. For $t=1$, $e_{ms}^{0}$ denotes the source-to-storage transfer initiated during the shift immediately preceding the scheduling horizon and is treated as known initial data. This timing convention prevents material dispatched from a source from being consumed by a production line before the required transfer delay has elapsed. 

The upper bound $B_s^{\max}$ represents the storage capacity of node $s$ and is selected from the maximum safe operating level specified by the facility for the corresponding silo or holding bin. Enforcing this bound prevents incoming mineral supply from exceeding the physical storage capacity of the facility, while the non-negativity of $B_s^t$ ensures that withdrawals by the production-line layer do not exceed the material available from existing inventory and completed source transfers.
The initial inventory at the beginning of the scheduling horizon, $t=0$, is specified by:
\begin{equation}
    B_s^0 = B_s^\mathrm{init}
    \quad \forall s \in \mathcal{S}
    \label{eq:initial_raw_storage}
\end{equation}
where $B_s^\mathrm{init}$ is the known initial inventory of raw material at storage node $s$ at the beginning of the planning horizon. Together, Equations \eqref{eq:storage_balance} and \eqref{eq:initial_raw_storage} connect the delayed mineral supply entering from the source layer to the material withdrawals required by the production-line layer while enforcing the storage and routing restrictions defined in Section \ref{sec:problem_definition}.

\subsubsection{Production Line Nodes}\label{subsubsec:line_constraints}
Production lines form the third layer of the network, where each line $\ell\in\mathcal{L}$ receives raw material from compatible storage nodes $s$ and converts it into finished product. Since material does not accumulate within this layer, the following material balance is enforced at each production line node $\ell$:
\begin{equation}
    \sum_{s \in \mathcal{S}_\ell} b^t_{s\ell}=\sum_{p \in \mathcal{P}_\ell}\sum_{w \in \mathcal{W}_p} x^t_{\ell pw} \quad \forall \ell \in \mathcal{L}, t \in \mathcal{T}
  \label{eq:line_bal}
\end{equation}
where the non-negative continuous variable $x^t_{\ell pw}$ represents the quantity of finished product $p$ produced on line $\ell$ directed to inventory location $w$ within shift $t$ after processing raw material input $b^t_{s\ell}$.

Line activation determines whether product $p$ is processed on line $\ell$ during shift $t$. This logical relationship is represented using GDP \cite{perez2023julia}. For each compatible product-line-shift combination, either product $p$ is active on line $\ell$ and may be produced, or it is inactive and its production quantity is zero. Let $Y_{\ell p}^t \in \{\mathrm{True},\mathrm{False}\}$ denote the Boolean activation variable, where $Y_{\ell p}^t=\mathrm{True}$ indicates that product $p$ is active on line $\ell$ during shift $t$. The corresponding GDP disjunction is:
\begin{equation}
    \Biggl[
    \begin{array}{c}
    Y_{\ell p}^t\\
    x^t_{\ell pw} \geq 0
    \end{array}
    \Biggr]
    \;\vee\;
    \Biggl[
    \begin{array}{c}
    \neg Y_{\ell p}^t\\
    x^t_{\ell pw} = 0
    \end{array}
    \Biggr]
    \quad
    \forall p\in\mathcal{P},\;w\in\mathcal{W}_p,\;\ell\in\mathcal{L}_p,\;t\in\mathcal{T}
    \label{eq:gdp_disjunction}
\end{equation}
The GDP disjunction is reformulated as a set of mixed-integer algebraic constraints by associating the Boolean variable $Y_{\ell p}^t$ with a binary variable $y_{\ell p}^t\in\{0,1\}$. When $y_{\ell p}^t=1$, the active disjunct is selected and product $p$ may be produced on line $\ell$ during shift $t$, subject to the remaining production constraints. Conversely, when $y_{\ell p}^t=0$, the inactive disjunct is selected and $x^t_{\ell pw}=0$ is enforced for every compatible inventory location $w\in\mathcal{W}_p$.

As discussed in Section \ref{sec:problem_definition}, transitioning between distinct product campaigns on the same line requires changeover activities that consume available operating time. Let $\delta_{\ell pp'}$ denote the changeover duration required on line $\ell$ when production transitions from product $p$ to product $p'$. Since the within-shift production sequence is not explicitly represented, the maximum pairwise changeover duration for each line is defined as:
\begin{equation}
    \delta^{{\max}}_{\ell} = \max_{\substack{(p,p')\in\mathcal{P}_\ell\times\mathcal{P}_\ell
    \\
    p\neq p'}} 
    \delta_{\ell pp'} \quad \forall \ell \in \mathcal{L}
    \label{eq:changeover_max}
\end{equation}
where $\delta^{\max}_{\ell}$ is the maximum changeover duration on line $\ell$. When multiple distinct products are processed on line $\ell$ during shift $t$, and each active product is processed in a single continuous run, the number of within-shift product transitions is one fewer than the number of active products. Since the production sequence is not explicitly modelled, the formulation assigns the maximum pairwise changeover duration on line $\delta_\ell^{\max}$ to each required transition, yielding the following upper bound on the total changeover time that may be required:
\begin{equation}
  \Delta^t_{\ell} \geq \delta^{\max}_{\ell} \left( \sum_{p \in \mathcal{P}_\ell} y^t_{\ell p} - 1 \right) \quad \forall \ell \in \mathcal{L}, t \in \mathcal{T}
  \label{eq:changeover_budget}
\end{equation}
where the non-negative continuous variable $\Delta_{\ell}^t$ is the total product-changeover time reserved for line $\ell$ during shift $t$. For any possible ordering of the active product campaigns, each pairwise transition duration satisfies $\delta_{\ell pp'}\leq\delta_{\ell}^{\max}$. Since the number of required transitions is one fewer than the number of active products, the right-hand side of Equation \eqref{eq:changeover_budget} provides an upper bound on the total changeover time required by any within-shift ordering. Requiring $\Delta_{\ell}^t$ to be no smaller than this bound is conservative because it may overestimate the time required by a known production sequence but cannot underestimate the time required by any possible within-shift product sequence.

To reflect practical industrial operating limitations, at most $\eta^{\max}$ products may be processed on a single production line within one shift:
\begin{equation}
  \sum_{p \in \mathcal{P}_\ell} y^t_{\ell p} \leq \eta^{\max}
  \quad \forall \ell \in \mathcal{L}, t \in \mathcal{T}
  \label{eq:max_products_shift}
\end{equation}
where $\eta^{\max}$ is the maximum number of distinct products permitted on a production line during a single shift, preventing an excessive number of product changeovers. This parameter should be selected based on historical operating practices and the number of campaigns that can realistically be completed within a shift after accounting for cleaning, purging, packaging changes, equipment configuration, and other transition requirements.

Although each time step $t$ provides a fixed window of operating time, not all time is available for production due to operator shift changes, equipment maintenance, and product changeovers. Let $D_{\mathrm{eff}}$ represent the effective production time available per shift. Production quantities are converted to processing times through the corresponding line and product-specific production rates, giving the line-capacity constraint:
\begin{equation}
    \sum_{p \in \mathcal{P}_\ell}\sum_{w \in \mathcal{W}_p} \frac{x^t_{\ell pw}}{r_{\ell p}} + \Delta_{\ell}^t \leq D_{\mathrm{eff}}
    \quad \forall \ell \in \mathcal{L}, t \in \mathcal{T}
    \label{eq:line_capacity}
\end{equation}
where $r_{\ell p}$ is the production rate for product $p$ on line $\ell$. Equation \eqref{eq:line_capacity} limits the productive and changeover time assigned to each line independently, but it does not capture competition between lines that share equipment. As described in Section \ref{sec:problem_definition}, some pairs of production lines share common equipment that can serve only one connected line at a time during active production. At the shift-level resolution, the aggregate productive time for which the paired lines occupy this resource is therefore constrained by:
\begin{equation}
    \sum_{p \in \mathcal{P}_{\ell_1}} \sum_{w \in \mathcal{W}_p} 
    \frac{x^{t}_{\ell_1pw}}{r_{\ell_1p}}
    +
    \sum_{p \in \mathcal{P}_{\ell_2}} \sum_{w \in \mathcal{W}_p} 
    \frac{x^{t}_{\ell_2pw}}{r_{\ell_2p}}
    \leq D_{\mathrm{eff}}
    \quad
    \forall (\ell_1,\ell_2)\in\mathcal{L}_{\mathrm{conflicts}},\; t\in\mathcal{T}
    \label{eq:conflict_constraint}
\end{equation}
where $\mathcal{L}_{\mathrm{conflicts}}$ is the set of unordered pairs of production lines that share a common resource, and each term $x^t_{\ell pw}/r_{\ell p}$ represents the portion of shift $t$ for which that resource is occupied by line $\ell$. Changeover-time variables are not included in Equation \eqref{eq:conflict_constraint} because the corresponding changeover activities are performed on line-specific equipment and do not occupy the shared resource. These activities are accounted for separately through $\Delta_{\ell}^t$ in the individual line-capacity constraint \eqref{eq:line_capacity}. Consequently, one line may undergo a changeover while the other continues to use the shared resource.

To prevent the model from creating a separate production campaign for a negligible quantity, a minimum productive duration is imposed whenever a product $p$ is activated on line $\ell$ during shift $t$:
\begin{equation}
  \sum_{w\in\mathcal{W}_p} \frac{x^t_{\ell pw}}{r_{\ell p}}
  \geq
  \tau_p^{\min} y^t_{\ell p}
  \quad
  \forall p \in \mathcal{P}, \ell \in \mathcal{L}_p, t \in \mathcal{T}
  \label{eq:min_run_time}
\end{equation}
where the left-hand side of Equation \eqref{eq:min_run_time} represents the productive operating time allocated to product $p$ on line $\ell$ during shift $t$. When $y_{\ell p}^t=1$, this time must be at least $\tau_p^{\min}$. The constraint prevents the optimizer from scheduling small residual quantities that would require a separate campaign but would not justify the associated line preparation, packaging configuration, product-flow stabilization, startup checks and subsequent transition to another product. The value of $\tau_p^{\min}$ should represent the shortest campaign that the facility would intentionally execute for product $p$, based on documented operating procedures or accepted operating heuristics. Changeover time is not included in $\tau_p^{\min}$ because it is accounted for separately through $\Delta_{\ell}^t$ in Equation \eqref{eq:changeover_budget}.

\subsubsection{Finished Product Inventory Nodes} \label{subsubsec:inventory_constraints}
Finished product inventory nodes $w\in\mathcal{W}$ represent the bulk-product holding areas, packaged-product staging locations, warehouse zones and other storage locations between the production lines and customer fulfillment layers. Once raw material is processed on line $\ell$, the resulting finished product $p$ is directed to a compatible inventory location $w\in\mathcal{W}_p$. Product-location compatibility may reflect mineral grade segregation, packaging type, contamination-prevention requirements, and established site-specific storage practices. 

For each compatible product-location pair, the inventory balance accounts for inventory carried forward from the preceding shift, finished product received from the production line layer and material shipped to customers:
\begin{equation}
  I^t_{pw} = I^{t-1}_{pw} + \sum_{\ell\in\mathcal{L}_p}x^t_{\ell pw} - \sum_{\substack{c\in\mathcal{C}\\(c,p)\in\mathcal{O}}}q^t_{wpc}
 \quad \forall p \in \mathcal{P}, w \in \mathcal{W}_p, t \in \mathcal{T}
  \label{eq:inv_balance}
\end{equation}
where $I^t_{pw}$ is the non-negative continuous inventory tracking variable representing the quantity of finished product $p$ at inventory location $w$ at the end of shift $t$, $x^t_{\ell pw}$ is the quantity of product $p$ produced on line $\ell$ allocated to location $w$ during shift $t$, and the non-negative continuous variable $q^t_{wpc}$ is the quantity of product $p$ shipped from inventory location $w$ to customer $c$ during shift $t$. Equation \eqref{eq:inv_balance} connects production decisions to customer fulfillment while ensuring that shipments cannot exceed the finished-product inventory available from preceding inventory and current-shift production. This balance assumes that product completed during a shift can become available for shipment within the same shift.

At the beginning of the scheduling horizon, the total available inventory of each product is specified as a known aggregate quantity that the model allocates among the compatible inventory locations:
\begin{equation}
    \sum_{w \in \mathcal{W}_p} I^{0}_{pw} = I^{\mathrm{tot}}_p \quad \forall p \in \mathcal{P}
    \label{eq:initial_inventory}
\end{equation}
where $I^{\mathrm{tot}}_p$ is the known total quantity of product $p$ available at the beginning of the horizon and $I^0_{pw}$ is the corresponding initial inventory variable allocated to inventory location $w$. 
This allocation does not represent a physical transfer of inventory at the beginning of the horizon; instead, it distributes the aggregate initial inventory among the compatible nodes used in the network representation.

Each inventory location is subject to lower and upper operating bounds:
\begin{equation}
  I^{min}_w \leq \sum_{p \in \mathcal{P}_w} I^t_{pw} \leq I^{max}_w \quad \forall w \in \mathcal{W}, t \in \mathcal{T}
  \label{eq:inventory_cap}
\end{equation}
where $I_w^{\min}$ and $I_w^{\max}$ are model parameters that denote the minimum and maximum allowable aggregate inventory levels at location $w$, respectively. Depending on the type of inventory node, these bounds may represent physical capacity, safety-stock requirements, or other operational inventory policies. Together, Equations \eqref{eq:inv_balance} to \eqref{eq:inventory_cap} connect finished production to customer shipments while enforcing product-location compatibility and the physical inventory limitations of a mineral-manufacturing facility.

\subsection{Customer Demand and Order Fulfillment Constraints}\label{subsec:demand_constraints}
The scheduling problem considers confirmed customer orders for which commercial terms, including pricing, have already been established outside the model. Consequently, customer requirements enter the formulation through specified product quantities and due dates rather than through price or revenue parameters. Let $\mathcal{O}\subseteq\mathcal{C}\times\mathcal{P}$ denote the set of customer-product pairs with existing confirmed orders. For each $(c,p)\in\mathcal{O}$, the non-negative parameter $d^t_{cp}$ denotes the quantity of product $p$ demanded by customer $c$ with a due date in shift $t$. 

To track the quantity of each order that has become eligible for fulfillment, the cumulative demand due up to and including shift $t$ is defined as:
\begin{equation}
D^t_{cp}=\sum_{i=1}^{t}d^i_{cp} 
\quad 
\forall (c,p)\in\mathcal{O},\;t\in\mathcal{T}
\label{eq:cumulative_demand}
\end{equation}
where the non-negative parameter $D^t_{cp}$ represents the cumulative demand due for customer-product pair $(c,p)$ through shift $t$. At the beginning of the scheduling horizon, cumulative demand arising within the horizon is initialized as: 
\begin{equation}
D^0_{cp}=0
\quad
\forall (c,p)\in\mathcal{O}
\label{eq:initial_cumulative_demand}
\end{equation}
Some customer-product pairs may enter the planning horizon with unmet demand from previous periods. Let the non-negative continuous variable $u^t_{cp}$ denote the cumulative unmet demand at the end of shift $t$ and $u^\mathrm{init}_{cp}$ denote the known initial backorder for customer-product pair $(c,p)$. Backorder is therefore defined as:
\begin{equation}
u^t_{cp} = u^\mathrm{init}_{cp}+D^t_{cp} - \sum_{i=1}^{t} \sum_{w\in\mathcal{W}_p} q^i_{wpc}
\quad
\forall (c,p)\in\mathcal{O},\;t\in\mathcal{T}
\label{eq:backorder}
\end{equation}
Equation \eqref{eq:backorder} accounts for both demand outstanding at the beginning of the horizon and demand that has become due during the horizon, factoring the cumulative quantity shipped to the customer. As $u_{cp}^t$ is constrained to be non-negative, cumulative shipments cannot exceed the initial backorder and demand due through the current shift. The formulation consequently prevents finished-product inventory from being allocated to future-dated demand before its corresponding due date.

\subsection{Lexicographic Optimization Approach}\label{subsec:lex_framework}
As established in Section \ref{sec:problem_definition}, the scheduling framework distinguishes between confirmed customer orders and forecast-based production targets. Confirmed orders represent accepted customer commitments, whereas forecast targets support anticipated future requirements. Commercial decisions concerning order acceptance and pricing occur outside the scope of this short-term scheduling problem. Consequently, the orders comprising $\mathcal{O}$ are not selected or prioritized according to their financial value. The adopted hierarchy instead reflects the facility's operating policy that fulfillment of confirmed orders must take precedence over production scheduled to improve forecast attainment.

The scheduling model employs this priority using sequential lexicographic optimization \cite{Ehrgott2005, Miettinen1998}. In lexicographic optimization, objectives are arranged in a strict order of importance and optimized through a sequence of individual solutions, each referred to as a phase. After a phase has been solved, constraints are added to preserve the performance achieved for its objective before the next, lower-priority objective is considered. Consequently, improvement in a lower-priority objective cannot be obtained by degrading performance associated with a previously optimized objective.

The framework consists of an initial customer-fulfillment phase followed by a sequence of forecast-oriented phases. To begin, Phase $1$ optimizes fulfillment of confirmed customer orders. Forecast-oriented Phase $n$ then considers the cumulative production target through forecast interval $n-1$, with one additional interval incorporated during each successive phase. This construction accommodates an arbitrary number $N$ of forecast intervals without changing the underlying formulation and optimization framework. Using $\succ$ to denote strict lexicographic priority, the resulting hierarchy is expressed as follows:
\begin{equation}
\text{Customer Fulfillment}
\succ
\text{Forecast interval }1
\succ
\text{Forecast intervals }1\text{--}2
\succ
\cdots
\succ
\text{Forecast intervals }1\text{--}N,
\label{eq:lexicographic_hierarchy}
\end{equation}

\subsubsection{Phase 1: Customer Fulfillment}
Phase $1$ establishes the highest priority customer service outcome by minimizing cumulative backorder over the scheduling horizon:
\begin{equation}
\min \quad
\sum_{(c,p)\in\mathcal{O}}
\sum_{t\in\mathcal{T}}
u^t_{cp},
\label{eq:obj_phase1}
\end{equation}
where $u^t_{cp}$ is the cumulative unmet demand for customer-product pair $(c,p)$ at the end of shift $t$. Backorder is minimized because it directly measures confirmed demand that has become due but remains unfulfilled. Since an unfulfilled quantity appears in $u^t_{cp}$ during every shift for which it remains outstanding, the objective penalizes both the quantity and duration of unmet demand. For example, one unit that remains unfulfilled for two shifts contributes to the objective in both shifts. Minimizing the sum therefore promotes earlier fulfillment of due orders without requiring time-dependent penalty weights.

Once Phase $1$ is solved, the resulting backorder levels are retained as upper bounds for every subsequent phase:
\begin{equation}
u^t_{cp}
\leq
\bar{u}^{t}_{cp},
\qquad
\forall (c,p)\in\mathcal{O}, t\in\mathcal{T}
\label{eq:phase1_backorder_pin}
\end{equation}
where $\bar{u}^{t}_{cp}$ denotes the cumulative backorder value obtained for customer-product pair $(c,p)$ at the end of shift $t$ in the optimal Phase $1$ solution. Each bound prevents a subsequent forecast-oriented phase from producing more backorder for any customer-product pair and shift than was obtained in Phase $1$. Collectively, these constraints preserve the customer-fulfillment profile attained throughout the horizon, rather than preserving only the aggregate value of the Phase $1$ objective. Forecast attainment can therefore be improved only through schedules that maintain the customer-service performance established in Phase $1$.

\subsubsection{Forecast-Oriented Phases}
After customer-fulfillment performance has been established, the forecast-oriented phases allocate the remaining feasible production capacity toward forecast-based production targets. As introduced in Section \ref{sec:problem_definition}, $\mathcal{K}=\{1,\ldots,N\}$ is the ordered set of forecast intervals, and $\mathcal{P}^{\mathrm{forecast}}$ is the set of products included in the forecast hierarchy. The target parameters are supplied externally by the demand-planning process and are treated as fixed inputs within the scheduling instance.

Forecast attainment is evaluated using cumulative scheduled production over the horizon rather than against customer shipments or ending inventory. Forecast-oriented Phase $n$, where $n\in \{2,\ldots,N+1\}$, considers forecast intervals $1,\ldots,n-1$. Phase $2$ considers the target for forecast interval $1$, Phase $3$ considers the cumulative target for forecast intervals $1$ and $2$, and each subsequent phase incorporates one additional interval. This ordering gives precedence to the earliest forecast requirements before progressively considering requirements over longer forecast intervals.

For each forecast-oriented Phase $n$, the effective cumulative production target for product $p$ is defined as:
\begin{equation}
    \widehat{\Phi}_p^{(n)} = \max\left\{0, \sum_{k=1}^{n-1}\Phi_p^k-\varphi_p\right\}
    \qquad
    \forall p\in\mathcal{P}^{\mathrm{forecast}},
    \quad
    n=2,\ldots,N+1
\label{eq:effective_forecast_target}
\end{equation}
where $\widehat{\Phi}_p^{(n)}$ is the effective cumulative production target for product $p$ in Phase $n$, $\Phi_p^k$ is the target production quantity for product $p$ in forecast interval $k$, and $\varphi_p$ is the production completed before the scheduling horizon within the current forecast cycle. Subtracting $\varphi_p$ prevents production completed before the scheduling horizon from being scheduled again, while the maximum operator sets the remaining target to zero when previously completed production already meets or exceeds the cumulative target. For each phase, a non-negative continuous variable $\psi^{(n)}_p$ represents the remaining production shortfall for product $p$ relative to its effective cumulative target:
\begin{equation}
    \psi^{(n)}_p
    \geq
    \widehat{\Phi}_p^{(n)} 
    -
    \sum_{w\in\mathcal{W}_p}\sum_{\ell\in\mathcal{L}_p}\sum_{t\in\mathcal{T}}x^t_{\ell pw},
    \qquad
    \psi^{(n)}_p\geq0,
    \quad
    \forall p\in\mathcal{P}^{\mathrm{forecast}},
    \;
    n=2,\ldots,N+1
    \label{eq:forecast_shortfall}
\end{equation}
The summation in Equation \eqref{eq:forecast_shortfall} represents the total quantity of product $p$ scheduled across all compatible production lines, finished-product inventory locations and shifts. Consequently, $\psi^{(n)}_p$ is positive only when scheduled production is insufficient to attain the effective target considered in Phase $n$. Phase $n$ therefore minimizes the aggregate production shortfall across all products included in the forecast:
\begin{equation}
\min \quad
\sum_{p\in\mathcal{P}^{\mathrm{forecast}}}
\psi^{(n)}_p
\qquad
\forall n=2,\ldots,N+1
\label{eq:obj_phase_n}
\end{equation}
The backorder bounds in Equation \eqref{eq:phase1_backorder_pin} remain active throughout every forecast-oriented phase. This ensures that forecast attainment can be improved without increasing the backorder of any customer-product pair during any shift. Beginning with the second forecast-oriented phase, the production achieved for each product in the preceding phase is also imposed as a lower bound:
\begin{equation}
\sum_{w\in\mathcal{W}_p}\sum_{\ell\in\mathcal{L}_p}\sum_{t\in\mathcal{T}}x^t_{\ell pw}
\geq
\bar{X}^{(n-1)}_p,
\qquad
\forall p\in\mathcal{P}^{\mathrm{forecast}},
\quad n=3,\ldots,N+1
\label{eq:forecast_phase_floor}
\end{equation}
where $\bar{X}^{(n-1)}_p$ is the production quantity achieved for product $p$ in Phase $n-1$. This product-level lower bound prevents a subsequent phase from improving attainment of a longer-term cumulative target by reallocating production away from a product supported in the preceding phase. The allocation of this preserved production across production lines, inventory locations, and shifts may change, provided that the total production achieved for each product does not decrease. Each successive phase therefore extends the forecast interval under consideration while preserving both the customer-service performance established in Phase $1$ and the product-level production achieved in the preceding forecast-oriented phase. With the lexicographic hierarchy established, the GDP disjunctions introduced in Section \ref{subsubsec:line_constraints} are reformulated such that the activation decisions are represented by binary variables and all remaining decision variables are continuous, yielding an MILP for each phase of the lexicographic solution procedure.

\section{Industrial Case Study Results} \label{sec:casestudy}
This section applies the sequential lexicographic MILP framework developed in Section \ref{sec:methodology} to operational data from a mineral-manufacturing facility in North America. The results evaluate computational performance, attainment of the prioritized scheduling objectives, and the behaviour induced by the principal operating constraints. To protect commercially sensitive information, facility-specific identifiers and operating details have been anonymized. All material quantities are reported in tonnes, preserving the operational relationships required for analysis without disclosing proprietary product or customer identities.

\subsection{Industrial Case Study and Operating Requirements} \label{subsec:case_description}
\begin{figure}[!htb]
    \centering
    \includegraphics[width=1\textwidth]{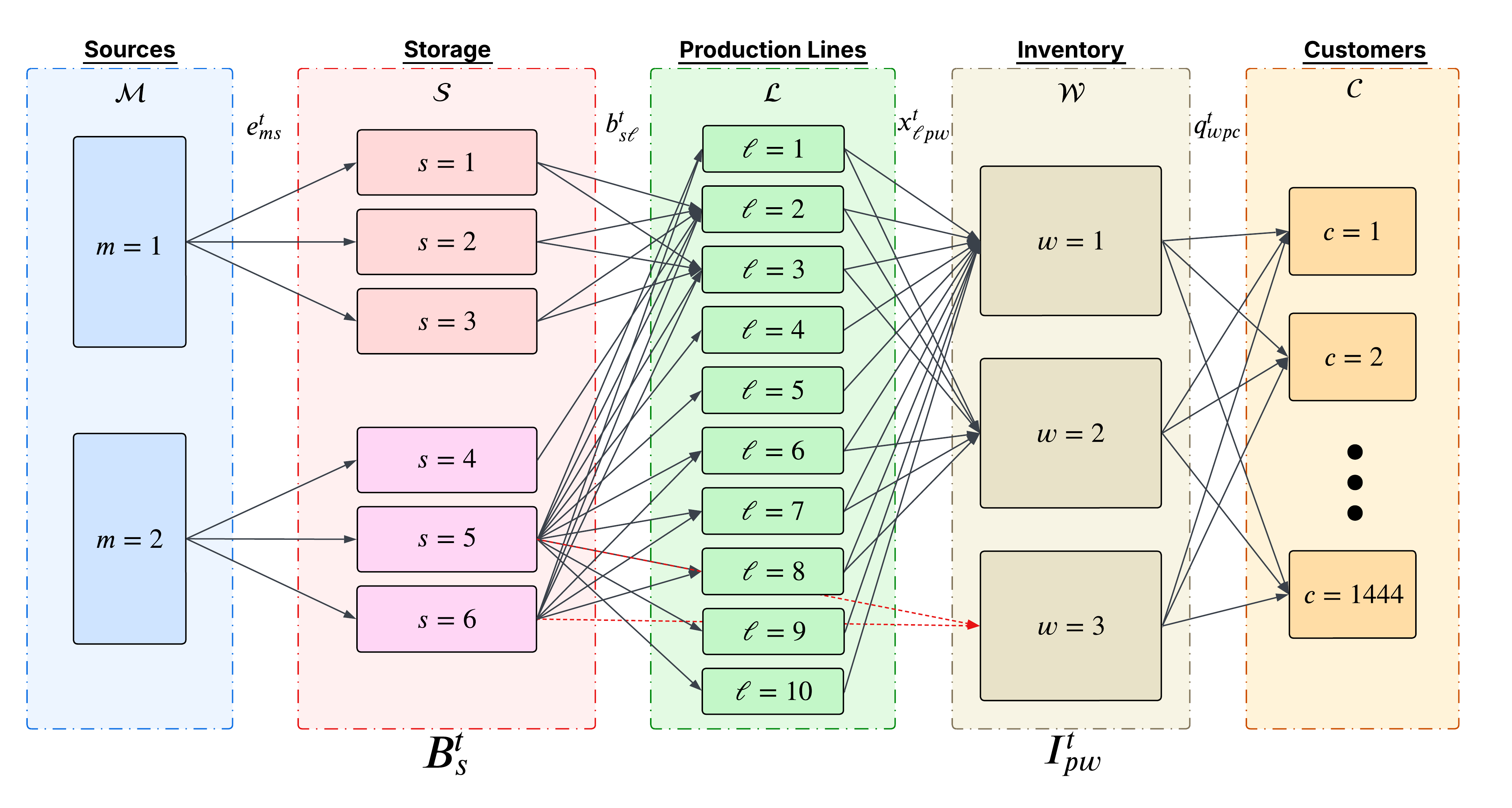}
    \caption{Case study production network illustrating the industrial supply chain flow.}
    \label{fig:PFD_case}
\end{figure}

The case study considers a multi-product mineral-manufacturing facility over a 28-day horizon comprising 84 eight-hour shifts. The production network contains two raw-material sources, six raw-material storage nodes, 10 production lines, three finished-product inventory locations, and 89 products. The principal dimensions of the production network and operational dataset are summarized in Table \ref{tab:case_study_sets}.
\begin{table}[!htb]
\centering
\caption{Cardinality of the principal sets for the industrial case study instance.}
\label{tab:case_study_sets}
\begin{tabular}{lll}
\hline
\textbf{Set} & \textbf{Description} & \textbf{Cardinality} \\
\hline
$\mathcal{M}$ & Raw-material sources & 2 \\
$\mathcal{S}$ & Raw-material storage nodes & 6 \\
$\mathcal{L}$ & Production lines & 10 \\
$\mathcal{W}$ & Finished-product inventory locations & 3 \\
$\mathcal{P}$ & Products & 89 \\
$\mathcal{C}$ & Customers & 1,444 \\
$\mathcal{O}$ & Customer--product demand pairs & 78 \\
$\mathcal{T}$ & Scheduling shifts & 84 \\
$\mathcal{K}$ & Forecast intervals & 2 \\
\hline
\end{tabular}
\end{table}
Figure \ref{fig:PFD_case} illustrates the case study network, which includes finite raw-material and finished-product storage capacities, product-dependent line compatibility and processing rates, and shared equipment between selected production lines.

The dataset contains 436 customer orders, which are aggregated by customer, product, and due shift into 78 active customer--product pairs. Two forecast intervals are considered, resulting in three optimization phases that prioritize customer fulfillment, current-period forecast attainment, and combined current- and subsequent-period forecast attainment. The case study also includes premium products that must be assigned to a designated inventory location and held for a specified period before newly produced material becomes available for shipment.

\subsubsection{Premium-Product Requirements}
Premium products $p\in\mathcal{P}^{\mathrm{prem}}$ must be assigned to the designated finished-product inventory location $w^{\mathrm{prem}}$:
\begin{equation}
\mathcal{W}_p = \{w^\mathrm{prem}\},
\quad
\forall p\in\mathcal{P}^{\mathrm{prem}}
\label{eq:prem_inv}
\end{equation}
The aggregate inventory of premium products at this location is limited by:
\begin{equation}
\sum_{p\in\mathcal{P}^{\mathrm{prem}}}
I^t_{p,w^{\mathrm{prem}}}
\leq
I^{\max}_{\mathrm{prem}}
\qquad
\forall t\in\mathcal{T}
\label{eq:cap_prem_indoor}
\end{equation}
Newly produced premium material must complete a hold of $H^{\mathrm{prem}}$ shifts before shipment, whereas initial inventory is assumed to have completed the required inspection before the scheduling horizon. The production shifts released by shift $t$ are:
\begin{equation}
\mathcal{R}^{\mathrm{prem}}(t)
=
\left\{
i\in\mathcal{T}
\;\middle|\;
i+H^{\mathrm{prem}}\leq t
\right\},
\quad
\forall t\in\mathcal{T}
\label{eq:premium_release_set}
\end{equation}
Cumulative premium-product shipments are therefore limited by the initial available inventory and production that has completed the required hold:
\begin{equation}
\sum_{i=1}^{t}
\sum_{\substack{c\in\mathcal{C}\\(c,p)\in\mathcal{O}}}
q^i_{w^{\mathrm{prem}}pc}
\leq
I^0_{p w^{\mathrm{prem}}}
+
\sum_{i\in\mathcal{R}^{\mathrm{prem}}(t)}
\sum_{\ell\in\mathcal{L}_p}
x^i_{\ell p w^{\mathrm{prem}}},
\qquad
\forall p\in\mathcal{P}^{\mathrm{prem}},
\;t\in\mathcal{T}
\label{eq:premium_buffer}
\end{equation}
This restriction prevents newly produced premium material from being shipped before the required hold has elapsed.

\subsection{Results}\label{subsec:results}
Following reformulation of the GDP disjunctions, the case study formulation results in a large-scale MILP containing 234,396 decision variables and 145,180 constraints. The model was implemented in \texttt{JuMP.jl} (v1.30.1), with the disjunctions expressed via \texttt{DisjunctiveProgramming.jl} (v0.6.1). The resulting MILP was solved using Gurobi 13.0.0 on a laptop equipped with a 3.30 GHz AMD Ryzen 7 PRO 8840HS processor and 64 GB of RAM running Windows 11. As summarized in Table \ref{tab:computational_performance}, the three optimization phases required 16.75 s of cumulative solver time, while the complete workflow required 18.53 s. This turnaround supports periodic re-optimization as updated production, inventory, and demand information becomes available.
\begin{table}[!htb]
\centering
\caption{Model size and computational performance for the industrial case study.}
\label{tab:computational_performance}
\begin{tabular}{lr}
\hline
\textbf{Metric} & \textbf{Value} \\
\hline
Continuous variables & 84,876 \\
Binary variables & 149,520 \\
Total variables & 234,396 \\
Constraints & 145,180 \\
\hline
Phase $1$ solve time (s) & 4.65 \\
Phase $2$ solve time (s) & 6.03 \\
Phase $3$ solve time (s) & 6.08 \\
Cumulative solver time (s) & 16.75 \\
Total workflow wall-clock time (s) & 18.53 \\
\hline
\end{tabular}
\end{table}

\begin{figure}[!htb]
    \centering
    \includegraphics[width=1\textwidth]{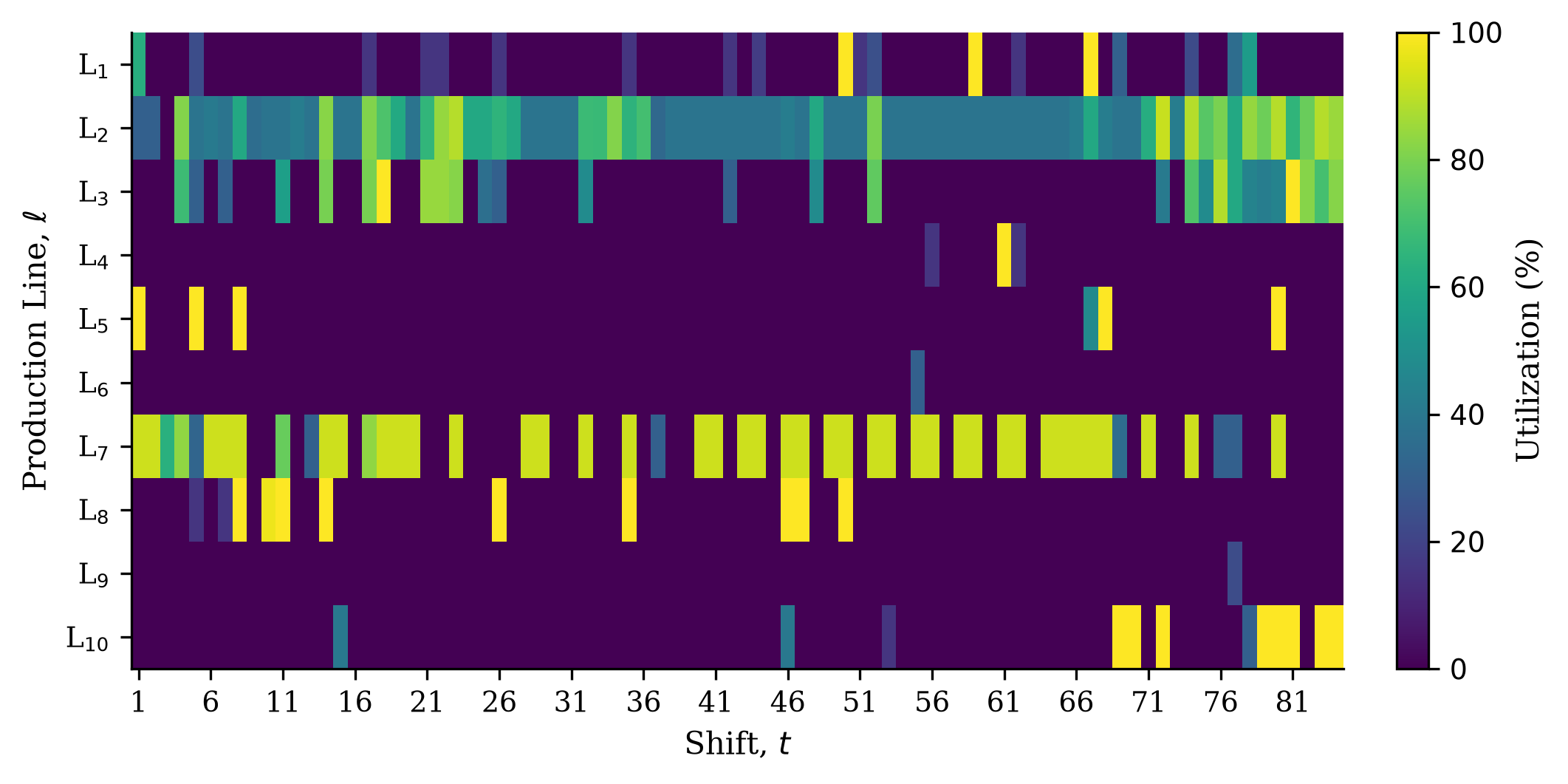}
    \caption{Production-line throughput utilization by shift over the 28-day scheduling horizon.}
    \label{fig:optimized_schedule}
\end{figure}

The final schedule produces 15,392 tonnes across 36 products and uses all 10 production lines, with 45 within-shift product changeovers. Lines $L_2$ and $L_3$ account for 79.4\% of total production, reflecting their compatible product portfolios and processing rates rather than an objective of distributing production uniformly. Figure \ref{fig:optimized_schedule} shows the resulting shift-level utilization of the production lines.

The schedule contains no incompatible product-line or product-inventory assignments, every active campaign satisfies its minimum duration, and no line exceeds the permitted number of products per shift. The combined productive time of each shared-resource line pair remains within the available shift duration. Premium-product output is assigned only to its designated inventory location, and shipments are restricted to inventory that has completed the required hold. These outcomes confirm that the principal operating constraints induce the intended scheduling behaviour.

\begin{figure}[!htb]
    \centering
    \includegraphics[width=0.85\textwidth]{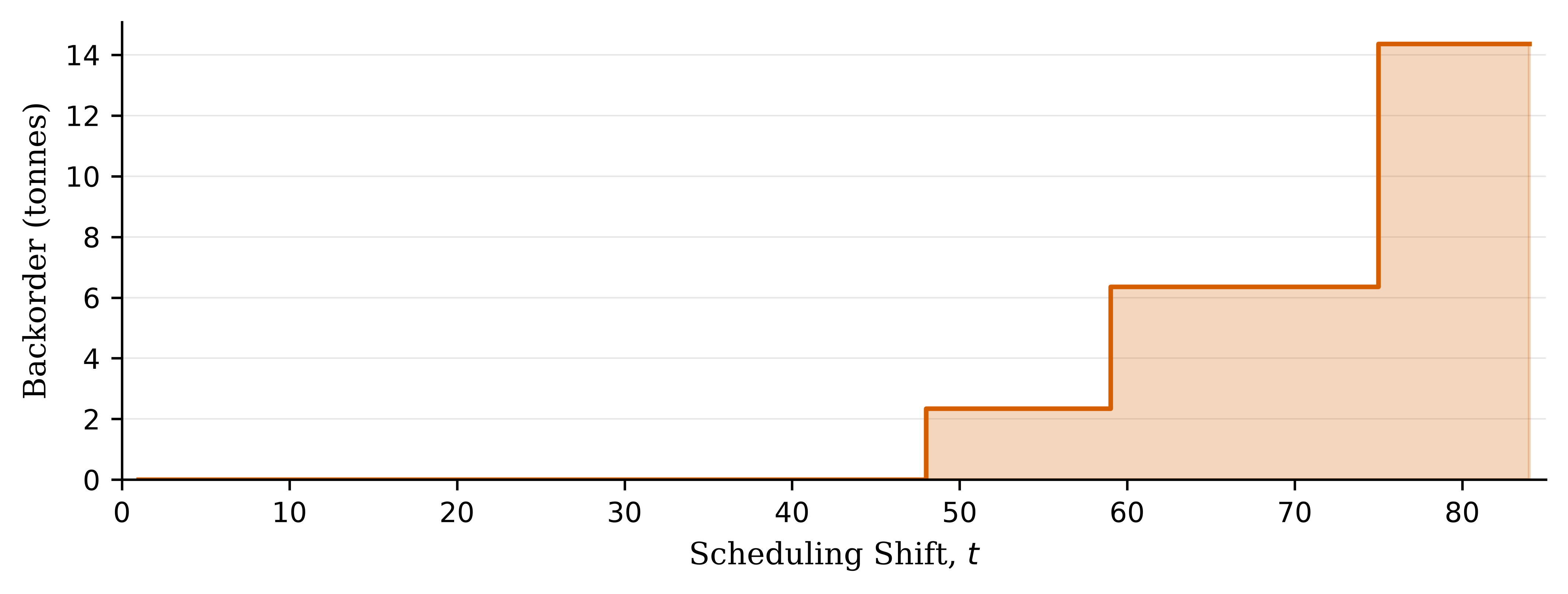}
    \caption{Aggregate customer backorder over the 84-shift scheduling horizon.}
    \label{fig:backorder_trajectory}
\end{figure}

Phase $1$ fulfills 8,657.38 of the 8,671.74 tonnes of confirmed demand, corresponding to 99.83\% fulfillment and a final backorder of 14.36 tonnes across four customer--product pairs. As shown in Figure \ref{fig:backorder_trajectory}, backorder remains zero through most of the horizon and arises only as the remaining requirements become due near its end. The same backorder profile is retained in Phases $2$ and $3$, confirming that forecast attainment is improved without degrading the higher-priority customer-service result.

\begin{figure}[!htb]
    \centering
    \includegraphics[width=0.9\textwidth]{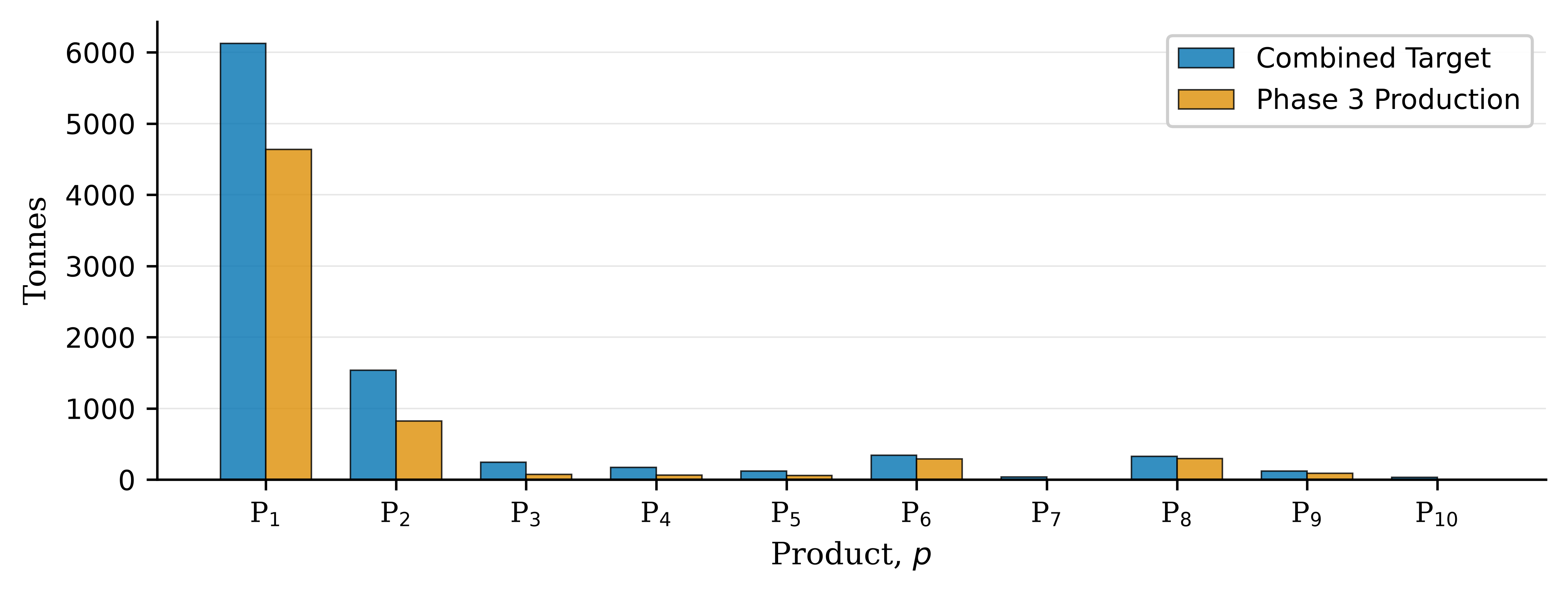}
    \caption{Combined forecast targets and scheduled production for the ten products contributing the largest Phase $3$ shortfalls.}
    \label{fig:phase3_forecast}
\end{figure}

After accounting for production completed before the scheduling horizon, the current-period target has an effective remaining requirement of 4,964.41 tonnes. Phase $2$ leaves a product-level shortfall of 55.22 tonnes, corresponding to 98.89\% attainment, with six products retaining positive shortfalls. Incorporating the subsequent-period forecast increases the effective cumulative target to 17,581.81 tonnes. Phase $3$ attains 83.68\% of this target, leaving a shortfall of 2,869.84 tonnes across 24 products. Table \ref{tab:forecast_performance} summarizes these results.
\begin{table}[!htb]
\centering
\caption{Forecast-target performance of the two forecast-oriented lexicographic phases.}
\label{tab:forecast_performance}
\begin{tabular}{lrr}
\hline
\textbf{Metric} & \textbf{Phase 2} & \textbf{Phase 3} \\
\hline
Current-period forecast target (tonnes) & 12,915.59 & 12,915.59 \\
Production-to-date (tonnes) & 8,476.02 & 8,476.02 \\
Remaining current-period requirement (tonnes) & 4,964.41 & -- \\
Subsequent-period forecast target (tonnes) & -- & 12,930.11 \\
Combined Phase $3$ target (tonnes) & -- & 17,581.81 \\
Forecast shortfall (tonnes) & 55.22 & 2,869.84 \\
Target attainment (\%) & 98.89 & 83.68 \\
Products in forecast set & 39 & 40 \\
Products with positive shortfall & 6 & 24 \\
\hline
\end{tabular}
\end{table}
Customer fulfillment remains at 99.83\% throughout both forecast-oriented phases, while the product-level production floors prevent Phase $3$ from reducing production established for the higher-priority current-period target. The three largest product shortfalls account for 82.5\% of the Phase $3$ deficit, as shown in Figure \ref{fig:phase3_forecast}. These results demonstrate the intended objective hierarchy: customer service is preserved first, followed by attainment of the current-period and subsequent-period forecast targets.

\begin{figure}[!htb]
    \centering
    \includegraphics[width=0.85\textwidth]{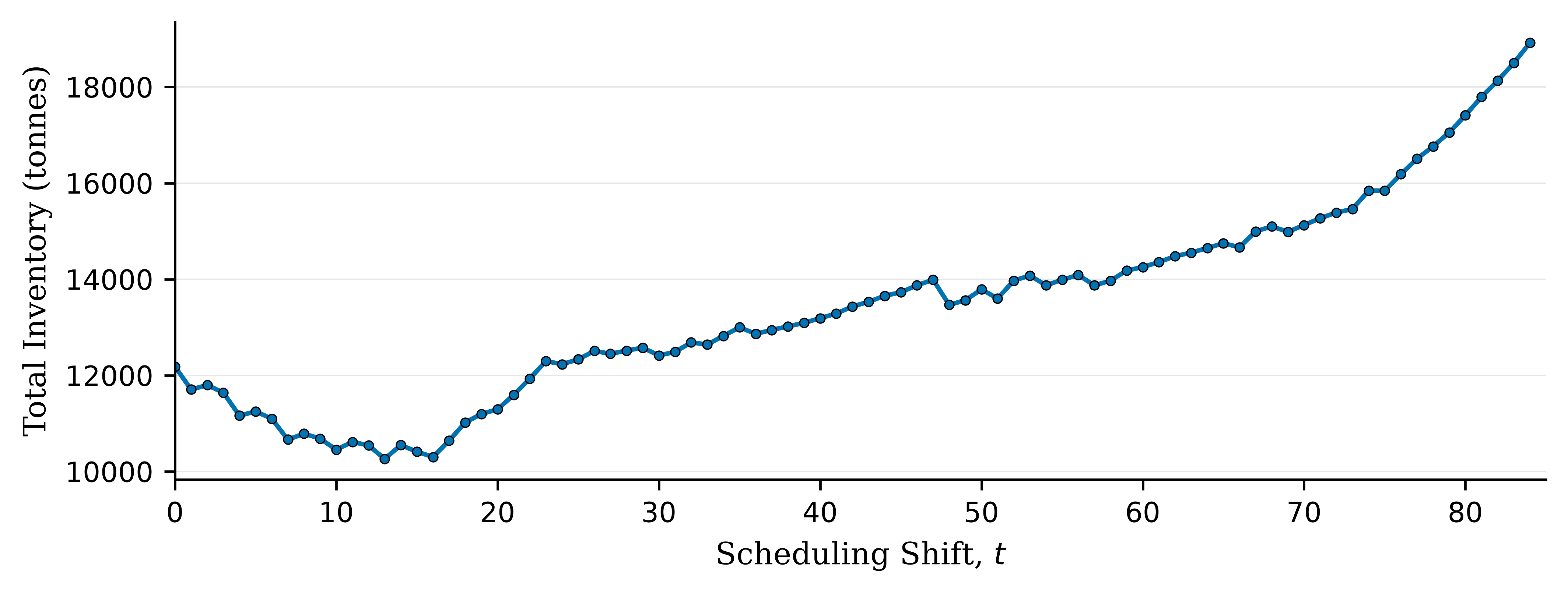}
    \caption{Aggregate finished-product inventory trajectory over the 84-shift scheduling horizon.}
    \label{fig:inventory_trajectory}
\end{figure}

Finished-product inventory decreases from 12,183 tonnes to a minimum of 10,262 tonnes at shift 13 as existing stock supports early customer requirements. Inventory subsequently increases to 18,918 tonnes by the end of the horizon as production directed toward the forecast targets exceeds shipments associated with the fixed customer-order set. Inventory remains within its operating bounds throughout the schedule, and the resulting trajectory is shown in Figure \ref{fig:inventory_trajectory}.

The late-horizon accumulation reflects the open-loop case study design: the optimization is solved once using the orders and plant state available at $t=0$, and no new orders are introduced during the 28-day horizon. The resulting trajectory is therefore conditional on the initial information and should not be interpreted as a recommended steady-state inventory target.

\section{Conclusions and Future Work} \label{sec:conclusion}
This work provides a shift-level production scheduling framework for multi-product mineral manufacturing. The formulation coordinates material flow, production, inventory, and customer fulfillment under detailed operating constraints. GDP represents logical line-product decisions, while a generalized lexicographic multi-objective strategy prioritizes customer fulfillment and sequentially addresses an arbitrary number of forecast intervals without degrading higher-priority objectives.

Applied to a 28-day industrial case study, the complete three-phase workflow required 18.53 s, fulfilled 99.83\% of customer demand, and attained 98.89\% and 83.68\% of the current-period and combined forecast requirements, respectively. These results demonstrate that detailed operating constraints and hierarchical production priorities can be addressed within a computationally tractable scheduling framework. Future work includes closed-loop re-optimization, higher-resolution sequencing, accounting for uncertainty, and economically driven objectives.

\section*{Acknowledgements}
We would like to acknowledge the financial support from both the Mitacs Accelerate Program (grant number IT47555) and the University of Waterloo.

\newpage

\section{Nomenclature}
\label{sec:notation}

The principal indices, sets, parameters, variables, and derived quantities used in the production scheduling formulation are summarized below. Material quantities are expressed in tonnes unless otherwise stated.

\subsection{Indices}

\begingroup
\small
\renewcommand{\arraystretch}{1.25}
\begin{longtable}{@{}p{0.14\textwidth}p{0.76\textwidth}@{}}
\caption{Indices used in the scheduling formulation.}
\label{tab:nomenclature_indices}\\
\toprule
\textbf{Index} & \textbf{Definition} \\
\midrule
\endfirsthead
\toprule
\textbf{Index} & \textbf{Definition} \\
\midrule
\endhead
$m$ & Raw-material source node \\
$s$ & Raw-material storage node \\
$\ell$ & Production line \\
$p,p'$ & Products; $p'$ is used as an alternative product index in changeover definitions \\
$w$ & Finished-product inventory location \\
$c$ & Customer \\
$t$ & Discrete scheduling period or production shift \\
$i$ & Auxiliary time index used in cumulative balances and release-set definitions \\
$k$ & Forecast-period index \\
$n$ & Lexicographic optimization-phase index, where Phase $n$ considers the first $n-1$ forecast intervals \\
\bottomrule
\end{longtable}
\endgroup

\subsection{Sets and Subsets}

\begingroup
\small
\renewcommand{\arraystretch}{1.25}
\begin{longtable}{@{}p{0.25\textwidth}p{0.65\textwidth}@{}}
\caption{Sets and subsets used in the scheduling formulation.}
\label{tab:nomenclature_sets}\\
\toprule
\textbf{Symbol} & \textbf{Definition} \\
\midrule
\endfirsthead
\toprule
\textbf{Symbol} & \textbf{Definition} \\
\midrule
\endhead
$\mathcal{M}$ & Set of raw-material source nodes. \\
$\mathcal{S}$ & Set of raw-material storage nodes. \\
$\mathcal{L}$ & Set of production lines. \\
$\mathcal{P}$ & Set of products. \\
$\mathcal{W}$ & Set of finished-product inventory locations. \\
$\mathcal{C}$ & Set of customers. \\
$\mathcal{T}$ & Ordered set of discrete scheduling periods or shifts. \\
$\mathcal{K}=\{1,\ldots,N\}$ & Ordered set of forecast intervals. \\
$\mathcal{O}\subseteq\mathcal{C}\times\mathcal{P}$ & Set of customer--product pairs having demand within the scheduling instance. \\
$\mathcal{M}_s\subseteq\mathcal{M}$ & Raw-material source nodes capable of supplying storage node $s$. \\
$\mathcal{S}_m\subseteq\mathcal{S}$ & Storage nodes that can receive raw material from source $m$. \\
$\mathcal{S}_{\ell}\subseteq\mathcal{S}$ & Storage nodes that can supply production line $\ell$. \\
$\mathcal{L}_s\subseteq\mathcal{L}$ & Production lines that can receive raw material from storage node $s$. \\
$\mathcal{L}_p\subseteq\mathcal{L}$ & Production lines compatible with product $p$. \\
$\mathcal{P}_{\ell}\subseteq\mathcal{P}$ & Products compatible with production line $\ell$. \\
$\mathcal{W}_p\subseteq\mathcal{W}$ & Finished-product inventory locations compatible with product $p$. \\
$\mathcal{P}_w\subseteq\mathcal{P}$ & Products compatible with finished-product inventory location $w$. \\
$\mathcal{P}^{\mathrm{forecast}}\subseteq\mathcal{P}$ & Products represented in the forecast-target hierarchy. \\
$\mathcal{P}^{\mathrm{prem}}\subseteq\mathcal{P}$ & Premium products subject to designated storage and release-hold requirements. \\
$\mathcal{L}_{\mathrm{conflicts}}\subseteq\mathcal{L}\times\mathcal{L}$ & Unordered pairs of production lines that share a common resource and are subject to a combined operating-time restriction. \\
$\mathcal{R}^{\mathrm{prem}}(t)$ & Production shifts whose premium-product output has completed the required hold and is available for shipment by shift $t$. \\
\bottomrule
\end{longtable}
\endgroup

\subsection{General and Temporal Parameters}

\begingroup
\small
\renewcommand{\arraystretch}{1.25}
\begin{longtable}{@{}p{0.25\textwidth}p{0.57\textwidth}p{0.11\textwidth}@{}}
\caption{General and temporal parameters.}
\label{tab:nomenclature_general_parameters}\\
\toprule
\textbf{Symbol} & \textbf{Definition} & \textbf{Unit} \\
\midrule
\endfirsthead
\toprule
\textbf{Symbol} & \textbf{Definition} & \textbf{Unit} \\
\midrule
\endhead
$N$ & Number of forecast intervals represented in the lexicographic hierarchy & -- \\
$D_{\mathrm{eff}}$ & Effective production time available on each production line during one shift & time/shift \\
$\eta^{\max}$ & Maximum number of distinct products that may be activated on one production line during one shift & -- \\
$H^{\mathrm{prem}}$ & Required premium-product hold duration, expressed in scheduling shifts & shifts \\
$w^{\mathrm{prem}}$ & Designated finished-product inventory location for premium products & -- \\
$E_m^t$ & Raw-material quantity supplied by source $m$ during shift $t$. & tonnes \\
\bottomrule
\end{longtable}
\endgroup

\subsection{Production and Changeover Parameters}

\begingroup
\small
\renewcommand{\arraystretch}{1.25}
\begin{longtable}{@{}p{0.25\textwidth}p{0.57\textwidth}p{0.11\textwidth}@{}}
\caption{Production and changeover parameters.}
\label{tab:nomenclature_production_parameters}\\
\toprule
\textbf{Symbol} & \textbf{Definition} & \textbf{Unit} \\
\midrule
\endfirsthead
\toprule
\textbf{Symbol} & \textbf{Definition} & \textbf{Unit} \\
\midrule
\endhead
$r_{\ell p}$ & Production rate of product $p$ on production line $\ell$, expressed using the same time basis as $D_{\mathrm{eff}}$ & tonnes/time \\
$\delta_{\ell pp'}$ & Changeover duration on line $\ell$ when transitioning from product $p$ to product $p'$ & time \\
$\delta_{\ell}^{\max}$ & Maximum pairwise changeover duration among products compatible with production line $\ell$ & time \\
$\tau_p^{\min}$ & Minimum production duration required whenever product $p$ is activated on a compatible line & time \\
\bottomrule
\end{longtable}
\endgroup

\subsection{Storage and Inventory Parameters}

\begingroup
\small
\renewcommand{\arraystretch}{1.25}
\begin{longtable}{@{}p{0.25\textwidth}p{0.57\textwidth}p{0.11\textwidth}@{}}
\caption{Raw-material storage and finished-product inventory parameters.}
\label{tab:nomenclature_inventory_parameters}\\
\toprule
\textbf{Symbol} & \textbf{Definition} & \textbf{Unit} \\
\midrule
\endfirsthead
\toprule
\textbf{Symbol} & \textbf{Definition} & \textbf{Unit} \\
\midrule
\endhead
$B_s^{\mathrm{init}}$ & Initial raw-material inventory at storage node $s$ & tonnes \\
$B_s^{\max}$ & Maximum raw-material inventory or physical capacity of storage node $s$ & tonnes \\
$I_p^{\mathrm{tot}}$ & Total initial finished-product inventory of product $p$ across its compatible inventory locations & tonnes \\
$I_w^{\min}$ & Minimum allowable aggregate finished-product inventory at location $w$ & tonnes \\
$I_w^{\max}$ & Maximum allowable aggregate finished-product inventory at location $w$ & tonnes \\
$I_{\mathrm{prem}}^{\max}$ & Maximum aggregate inventory of premium products at the designated location $w^{\mathrm{prem}}$ & tonnes \\
\bottomrule
\end{longtable}
\endgroup

\subsection{Demand and Forecast Parameters}

\begingroup
\small
\renewcommand{\arraystretch}{1.25}
\begin{longtable}{@{}p{0.25\textwidth}p{0.57\textwidth}p{0.11\textwidth}@{}}
\caption{Customer-demand and forecast parameters.}
\label{tab:nomenclature_demand_parameters}\\
\toprule
\textbf{Symbol} & \textbf{Definition} & \textbf{Unit} \\
\midrule
\endfirsthead
\toprule
\textbf{Symbol} & \textbf{Definition} & \textbf{Unit} \\
\midrule
\endhead
$d_{cp}^t$ & Demand for product $p$ from customer $c$ that becomes due during shift $t$ & tonnes \\
$D_{cp}^t$ & Cumulative demand for customer--product pair $(c,p)$ due through shift $t$ & tonnes \\
$u_{cp}^{\mathrm{init}}$ & Initial unmet demand or backorder for customer--product pair $(c,p)$ at the beginning of the scheduling horizon & tonnes \\
$\Phi_p^k$ & Production target for product $p$ in forecast interval $k$ & tonnes \\
$\Phi_p^1$ & Current-period production target for product $p$ in case study implementation & tonnes \\
$\Phi_p^2$ & Subsequent-period production target for product $p$ in the case study implementation & tonnes \\
$\varphi_p$ & Production of product $p$ completed before the scheduling horizon but within the current forecast interval & tonnes \\
$\widehat{\Phi}_p^{(n)}$ & Effective cumulative production target for product $p$ in Phase $n$ & tonnes \\
\bottomrule
\end{longtable}
\endgroup

\subsection{Continuous Decision Variables}

\begingroup
\small
\renewcommand{\arraystretch}{1.25}
\begin{longtable}{@{}p{0.25\textwidth}p{0.57\textwidth}p{0.11\textwidth}@{}}
\caption{Continuous decision variables.}
\label{tab:nomenclature_continuous_variables}\\
\toprule
\textbf{Symbol} & \textbf{Definition and domain} & \textbf{Unit} \\
\midrule
\endfirsthead
\toprule
\textbf{Symbol} & \textbf{Definition and domain} & \textbf{Unit} \\
\midrule
\endhead
$e_{ms}^t\geq0$ & Raw-material quantity transferred from source $m$ to compatible storage node $s$ during shift $t$ & tonnes \\
$B_s^t\in[0,B_s^{\max}]$ & Raw-material inventory at storage node $s$ at the end of shift $t$ & tonnes \\
$b^t_{s\ell}\geq0$ & Raw-material quantity transferred from storage node $s$ to compatible production line $\ell$ during shift $t$ & tonnes \\
$x^t_{\ell pw}\geq0$ & Quantity of product $p$ produced on line $\ell$ and allocated to compatible finished-product inventory location $w$ during shift $t$ & tonnes \\
$\Delta_{\ell}^t\geq0$ & Changeover time reserved on production line $\ell$ during shift $t$ & time \\
$I^t_{pw}\geq0$ & Inventory of product $p$ at compatible finished-product inventory location $w$ at the end of shift $t$ & tonnes \\
$I_{pw}^{0}\geq0$ & Initial inventory allocation of product $p$ to compatible inventory location $w$ & tonnes \\
$q^t_{wpc}\geq0$ & Quantity of product $p$ shipped from inventory location $w$ to customer $c$ during shift $t$ & tonnes \\
$u_{cp}^t\geq0$ & Cumulative unmet demand or backorder for customer--product pair $(c,p)$ at the end of shift $t$ & tonnes \\
$\psi_p^{(n)}\geq0$ & Product-level forecast shortfall minimized during lexicographic Phase $n$ & tonnes \\
\bottomrule
\end{longtable}
\endgroup

\subsection{Logical Variables and Reformulation Symbols}

\begingroup
\small
\renewcommand{\arraystretch}{1.25}
\begin{longtable}{@{}p{0.25\textwidth}p{0.57\textwidth}p{0.11\textwidth}@{}}
\caption{Logical activation variables and their binary reformulation.}
\label{tab:nomenclature_logical_variables}\\
\toprule
\textbf{Symbol} & \textbf{Definition and domain} & \textbf{Unit} \\
\midrule
\endfirsthead
\toprule
\textbf{Symbol} & \textbf{Definition and domain} & \textbf{Unit} \\
\midrule
\endhead
$Y_{\ell p}^t\in\{\mathrm{True},\mathrm{False}\}$ & Boolean GDP variable indicating whether product $p$ is activated on production line $\ell$ during shift $t$ & -- \\
$y_{\ell p}^t\in\{0,1\}$ & Binary reformulation of $Y_{\ell p}^t$; equal to one when product $p$ is active on line $\ell$ during shift $t$ and zero otherwise & -- \\
\bottomrule
\end{longtable}
\endgroup

\subsection{Retained Phase-Solution Quantities and Operators}

\begingroup
\small
\renewcommand{\arraystretch}{1.25}
\begin{longtable}{@{}p{0.25\textwidth}p{0.65\textwidth}@{}}
\caption{Quantities retained between lexicographic phases and mathematical operators.}
\label{tab:nomenclature_phase_quantities}\\
\toprule
\textbf{Symbol} & \textbf{Definition} \\
\midrule
\endfirsthead
\toprule
\textbf{Symbol} & \textbf{Definition} \\
\midrule
\endhead
$\bar{u}_{cp}^t$ & Backorder for customer--product pair $(c,p)$ at shift $t$ obtained from the optimal Phase $1$ solution and imposed as an upper bound in subsequent phases. \\
$\bar{X}_p^{(n-1)}$ & Total production of product $p$ obtained in forecast-oriented Phase $n-1$ and imposed as a product-level lower bound in Phase $n$. \\
$\succ$ & Strict lexicographic-priority relation. \\
$|\mathcal{A}|$ & Cardinality of a generic set $\mathcal{A}$. \\
\bottomrule
\end{longtable}
\endgroup

\newpage
\bibliography{references}

\end{document}